\documentclass[final]{proc-l} 
\begin{document}


\newcommand{\LL}{{\mathbb{L}}}
\newcommand{\II}{{\mathbb{I}}}
\newcommand{\JJ}{{\mathbb{J}}}

\newcommand{\calP}{{\mathcal P}}
\newcommand{\calM}{{\mathcal M}}
\newcommand{\calN}{{\mathcal N}}
\newcommand{\calQ}{{\cal Q}}

\newcommand{\sX}{{\mathbf{X}}}
\newcommand{\sF}{{\Hat{\mathbf{X}}}}

\newcommand{\sR}{{\mathbb{R}}}
\newcommand{\NN}{{\mathbb{N}}}

\newcommand{\sC}{{\mathbf{C}}}
\newcommand{\bx}{{\bf{x}}}
\newcommand{\by}{{\bf{y}}}

\newcommand{\lbl}[1]{\label{#1}}
\newcommand{\rf}[1]{(\ref{#1})}
\newcommand{\LVF}{\index{Varadhan Functional}Varadhan Functional}
\newcommand{\const}{\mathit{const}}
\newcommand{\La}{\Lambda}
\newcommand{\la}{\lambda}
\newcommand{\eps}{\epsilon}

\theoremstyle{definition}
\newtheorem{Definition}{Definition}[section]
\newtheorem{Claim}{Claim}[section]
\newtheorem{Exercise}{Exercise}[section]
\newtheorem{Condition}{Condition}[section]
\newtheorem{Problem}{Problem}
\theoremstyle{remark}
\newtheorem{Remark}{Remark}[section]
\theoremstyle{plain}
\newtheorem{Theorem}{Theorem}[section] 
\newtheorem{Corollary}[Theorem]{Corollary}
\newtheorem{Lemma}[Theorem]{Lemma}
\newtheorem{Example}{Example}[section]
\newtheorem{Proposition}[Theorem]{Proposition}

\newenvironment{proofof}[1]{\noindent {\bf Proof of #1.\/}}{\qed\vskip 0.1in}

\author{Harold Bell}
\address{
 Department of Mathematics \\
University of Cincinnati\\
PO Box 210025\\
Cincinnati, OH 45221--0025\\
}
\email{bellh@math.uc.edu}
\author{ Wlodzimierz  Bryc }
\address{ Department of Mathematics \\
University of Cincinnati\\
PO Box 210025\\
Cincinnati, OH 45221--0025}
\email{brycwz@email.uc.edu}
 \urladdr{http://http://ucms02.csm.uc.edu/preprint/ld-abs}

\title{Variational representations of Varadhan Functionals}
\date{June  2, 1999. Revised: November 1, 1999}

\begin{abstract} 
Motivated by the theory of large deviations, 
we introduce a class of non-negative non-linear functionals that 
have a variational ``rate function" representation.
\end{abstract}
\maketitle
\section{Introduction}

Let $(\sX,d)$ be a Polish space with metric $d()$ and let $\sC_b(\sX)$
denote the space of all bounded continuous functions $F:\sX\to\sR$.
In his work on large deviations of probability measures $\mu_n$,
 Varadhan \cite{Varadhan(1966)} introduced a class of non-linear functionals $\LL$ 
defined by 
\begin{equation} \label{LDP} 
\LL(F)=\lim_{n\to\infty} \frac1n\log\int_{\sX}\exp (n F(\bx)) d\mu_n
\end{equation}
 and used  the
large deviations principle of $\mu_n$
to prove the  variational
representation
\begin{equation} \label{(1.7)}   \LL (F)=L_0+\sup_{\bx\in \sX} \{ F(\bx)-\II(\bx)\},
\end{equation}  
where $\II:\sX\to[0,\infty]$ is the {\em rate function} governing the large deviations,
 and $L_0:=\LL(0)=0$.

Several authors 
\cite{Akian99,Bryc90,deAcosta(1985),OBrien(1996),OBrien-Vervaat(1995),Puhalskii(1997)}
 abstracted
non-probabilistic components from the theory of large deviations.
In particular,
in \cite{Bryc90}, see also \cite[Theorem 3.1]{OBrien-Vervaat(1995)} we give  
 conditions which imply  the rate function representation (\ref{(1.7)})
 when the limit (\ref{LDP}) exists,
and we show that the rate function is determined from the dual formula
\begin{equation}\label{(1.6)}   \II(\bx)=\LL(0)+\sup_{F\in \sC_b(\sX)}\{ F(\bx)-\LL (F)\}.
\end{equation}
In fact, one can reverse Varadhan's approach, and 
show that large deviations of probability measures 
$\mu_n$ follow from the variational representation (\ref{(1.7)}) for (\ref{LDP}), 
see \cite[Theorem 1.2.3]{Dupuis-Ellis(1997)}. In this context we have $\mu_n(\sX)=1$ which implies
$\LL(0)=0$ in (\ref{(1.6)}) and correspondingly $L_0=0$ in (\ref{(1.7)}).

 ``Asymptotic values" in \cite{Bryc90} are essentially what we call
\LVF s here; the theorems in that paper are not entirely satisfying because
the  assumptions 
are in terms of the underlying probability measures.   
In this paper we present a more satisfying 
approach which relies on the theory of probability for motivation purposes
only. 
\begin{Definition}
A function $\LL:\sC_b(\sX)\to\sR$ is  a \LVF\ 
if the following  conditions are satisfied.
\begin{eqnarray}
\label{V1.2} \hbox{If $F\leq G$ then $\LL(F)\leq \LL(G)$ for all } F, G\in\sC_b(\sX)&&\\
\label{V1.3}\LL(F+\const)=\LL(F)+\const \hbox{ for all }F\in\sC_b(\sX),\const\in\sR&&
\end{eqnarray}
\end{Definition}
Expression (\ref{LDP}) provides an example of \LVF, if the limit exists.
Another example is given by variational representation \rf{(1.7)}.

Condition \rf{V1.2} is  equivalent to  $\LL(F\vee G) \geq \LL(F)\vee \LL(G)$,
where $a\vee b$ denotes the maximum of two numbers. \LVF s like \rf{LDP} satisfy a stronger
condition.
\begin{Definition}
A \LVF\ $\LL$ is {\em maximal} \index{\LVF*maximal} if 
 $\LL(\cdot)$ is a lattice homomorphism
\begin{equation} \label{MAX}
 \LL(F\vee G)= \LL(F)\vee \LL(G).
\end{equation}
\end{Definition}
It is easy to see that each \LVF\ $\LL(\cdot)$ satisfies the Lipschitz condition  
$|\LL(F)-\LL(G)|\leq \|F-G\|_\infty$, compare \rf{(2.2)}.  Thus $\LL$ is a continuous
mapping from the Banach space $\sC_b(\sX)$ of all bounded continuous functions into the real
line. 
We will need the following stronger
continuity assumption, motivated by the definition of the
countable additivity of measures.
\begin{Definition}
A \LVF\ is {\em $\sigma$-continuous} if the following
condition is satisfied.
\begin{equation}
\label{V1.4} \hbox{ If $F_n\searrow 0$ then } \LL(F_n)\to\LL(0).
\end{equation}
\end{Definition}
Notice that if $\sX$ is compact, then  by Dini's theorem and the
Lipschitz property, all \LVF s are $\sigma$-continuous.

Maximal \LVF s are convex; this follows from the proof of Theorem \ref{T.1.2}, which
 shows that  formula \rf{(1.7)} holds
true for all \LVF s when the supremum is extended to all $\bx$ in the \v{C}ech-Stone
compactification of $\sX$.

A simple example of convex and
maximal but not $\sigma$-continuous \LVF\ is 
$\LL(F)=\limsup_{x\to\infty}F(x)$, where  $F\in\sC_b(\sR)$. This \LVF\ cannot be represented
by variational formula \rf{(1.7)}. Indeed,
\rf{(1.7)} implies that $\II(\bx)\geq F(\bx)-\LL(F)=F(\bx)$ for all
$F\in\sC_b(\sR)$
 that vanish at $\infty$;  hence $\II(\bx)=\infty$ for all
$\bx\in\sR$ and \rf{(1.7)}  gives $\LL(F)=-\infty$ for all $F\in\sC_b(\sR)$, a contradiction. 

An example of a convex  and
$\sigma$-continuous but not maximal \LVF\ is 
$\LL(F)=\log\int_{\sX} \exp F(\bx) \nu(d\bx)$, where $\nu$ is a finite non-negative measure.

\section{Variational
representations}

The main result of this paper is the following.
\begin{Theorem}\label{T.1.2} If a
maximal \LVF\ $\LL:\sC_b(\sX)\to\sR$ is 
$\sigma$-continuous,
then there is $L_0\in\sR$ such that
 variational representation \rf{(1.7)} holds true and
the rate function  $\II:\sX\to[0,\infty]$  is given by the dual formula \rf{(1.6)}. 
Furthermore,  $\II(\cdot)$ is a tight rate function: sets $\II^{-1}([0,a])\subset\sX$ are compact
for all $a>0$.
\end{Theorem}

The next result is closely related to Bryc \cite[Theorem T.1.1]{Bryc90} and
Deuschel \& Stroock  \cite[Theorem 5.1.6]{D-S}.
 Denote by $\calP(\sX)$  the metric space (with Prokhorov metric)
of all
probability measures on  a Polish space $\sX$ with the
Borel $\sigma$-field generated by all open sets.
\begin{Theorem}\label{mini-T.1.1}
 If a
convex \LVF\ $\LL:\sC_b(\sX)\to\sR$ is 
$\sigma$-continuous, then
there is a lower semicontinuous function $\JJ:\calP(\sX)\to [0,\infty]$ 
and a constant $L_0$ such that
such that  
\begin{equation}\label{CVX}
\LL(F)=L_0+\sup_{\mu\in\calP}\{ \int F d\mu -\JJ(\mu)\}
\end{equation}
for all bounded continuous functions $F$. 
\end{Theorem}

A well known example in large deviations is
 the convex $\sigma$-continuous functional
$\LL(F):=\log\int\exp F(\bx)\nu(d\bx)$ with the rate function in \rf{CVX} given by the
relative entropy
functional
$$\JJ(\mu)=\left\{\begin{array}{cl}
\int \log \frac{d\mu}{d\nu} d\mu & \hbox{ if $\mu\ll\nu$ is absolutely continuous} \\ 
\infty &\hbox{otherwise.}\end{array}\right.$$ 

\begin{Remark}
Deuschel \& Stroock \cite[Section 5.1]{D-S}  consider 
 convex functionals $\Phi:\sC_b(\sX)\to\sR$ such that
 $\Phi(\const)=\const$.
Such functionals satisfy condition \rf{V1.3}. Indeed,
write $F+\const$ as a convex combination
$$F+\const=(1-\theta)F+\frac\theta2\left(\frac{2\const}{\theta}\right)
+\frac{\theta}{2}(2F),
$$
where $0<\theta<1$.
Using convexity and 
 $\Phi(\const)=\const$ we get
$\Phi(F+\const)\leq \Phi(F)+\const+\theta(\frac{\Phi(2F)}{2}-\Phi(F))$.
Since $\theta>0$ is arbitrary this proves that $\Phi(F+\const)\leq \Phi(F)+\const$.
By routine symmetry considerations 
(replacing $F\mapsto F-\const$, and then   $const\mapsto -\const$), \rf{V1.3} follows.
\end{Remark}

\section{Proofs}
Let $L_0:=\LL(0)$. Passing to $\LL'(F):=\LL(F)-L_0$ if necessary, without losing generality we assume
$\LL(0)=0$.

\begin{Lemma}\label{Cech} Let $\sF$ be a compact Hausdorff space. Suppose $\sX\subset\sF$ is a
separable metric space in the relative topology. If $\bx_0\in\sF\setminus \sX$ then
there are bounded continuous functions $F_n:\sF\to\sR$ such that
\begin{itemize}
\item[(i)] $F_n(\bx)\searrow 0$ for all $\bx\in\sX$.
\item[(ii)] $F_n(\bx_0)=1$ for all $n\in\NN$.
\end{itemize}
\end{Lemma}
\begin{proof}
Since $\sF$ is Hausdorff, for every $\bx\in\sX$ there is an open set $U_\bx\ni\bx$ such that
its closure $\bar{U}_\bx$ does not contain $\bx_0$.

By Lindel\"of property for separable metric space $\sX$,
 there is a countable subcover $\{U_n\}$ of  $\{U_{\bx}\}$. 

A compact Hausdorff space $\sF$ is normal. So there are continuous functions
$\phi_n:\sF\to \sR$ such that $\phi_n\left|_{\bar{U}_n}\right.=0$ and $\phi_n(\bx_0)=1$.

To end the proof take $F_n(\bx)=\min_{1\leq k\leq n}\phi_k(\bx)$.
 
\end{proof}
The following lemma is contained implicitly in \cite[Theorem T.1.2]{Bryc90}.
\begin{Lemma}\label{L-T.1.2} Theorem \ref{T.1.2} holds true for compact $\sX$.
\end{Lemma}
\begin{proof} 

Let $\II(\cdot)$ be defined by \rf{(1.6)}. Thus $\II(\bx)\geq F(\bx)-\LL(F)$
which implies
$   \LL(F)\geq \sup_{\bx\in \sX} \{ F(\bx)-\II(\bx)\}$.
To end the proof  we need therefore to establish the converse inequality.
Fix a bounded continuous function 
 $F\in \sC_b(\sX)$ and $\eps>0$. Let 
$s=\sup_{\bx\in \sX} 
\{ F(\bx)-\II(\bx)\}$. Clearly $
F(\bx)-\II(\bx)\leq s\leq \LL(F)$. By  \rf{(1.6)} again, for
every $\bx\in \sX$, there is $F_\bx\in
\sC_b(\sX)$ such 
that $\II(\bx)<F_\bx(\bx)-\LL(F_\bx)+\eps$. Therefore
   $$F(\bx)\leq s+\II(\bx)<s+\eps+F_\bx(\bx)-\LL(F_\bx)$$
This means that the sets $U_\bx=\{ \by\in \sX: F(\by)-F_\bx(\by)<s+\eps-\LL(F_\bx)\} $ 
 form an open covering of $\sX$. Using compactness of 
$\sX$, we choose a finite covering $U_{\bx(1)}, \dots , U_{\bx(k)}$.
Then, writing $F_i=F_{\bx(i)}$ we have 
   $$F(\bx)< \max_{1 \leq i \leq k}\{F_{i}(\bx)-\LL(F_{i})\}+s+\eps
$$
for all  $\bx\in \sX$.

Using  \rf{V1.2}, \rf{V1.3}, and  \rf{MAX} we have
$$\LL(F) \leq \LL\left( \max_{1 \leq i \leq k}
\left\{F_{i}-\LL(F_{i})\right\} + s+\eps\right)=
\LL\left(\max_{ i}\left\{F_{i}-\LL(F_{i})\right\}\right) + s+\eps=$$
$$\max_{i}\left\{\LL\left(F_{i}-\LL(F_{i})\right)\right\} + s+\eps$$
Since \rf{V1.3} implies $\LL\left(F_{i}-\LL(F_{i})\right)=\LL\left(F_{i}\right)-\LL(F_{i})=0$
this shows that $s\leq \LL(F)<s+\eps$. Therefore $\LL(F)=s$, 
proving \rf{(1.7)}.
\end{proof} 

\begin{proofof}{Theorem \protect{\ref{T.1.2}}}
Let $\sF$ be  the \v{C}ech-Stone compactification of $\sX$. 
Since the inclusion $\sX\subset\sF$ is continuous, 
we define $\Hat{\LL}: \sC_b(\sF)\to\sR$ by $\Hat{\LL}(\Hat{F}):=\LL(\Hat{F}|_\sX)$.
It is clear that $\Hat{\LL}$ is a maximal \LVF, so
by Lemma \ref{L-T.1.2} there is $\II:\sF\to[0,\infty]$ such that
$\Hat{\LL}(\Hat{F})=\sup\{\Hat{F}(\bx)-\II(\bx): \bx\in\sF\}$.

Using $\sigma$-continuity\rf{V1.4}  it is  easy to check
that $\II(\bx)=\infty$ for all
$\bx\in\sF\setminus\sX$.
Indeed, given $\bx_0\in\sF\setminus\sX$ by Lemma \ref{Cech} there are
$F_n\in\sC_b(\sF)$  such that 
$F_n\searrow 0$ on $\sX$, but $F_n(\bx_0)=C>0$.
Then from \rf{(1.6)} we get $\II(\bx_0)\geq
\Hat{\LL}(0)+F_n(\bx_0)-\Hat{\LL}(F_n)\to \Hat{\LL}(0)+C$.  Since $C>0$ is arbitrary,
$\II(\bx_0)=\infty$.

This shows that $\Hat{\LL}(\Hat{F})=\sup\{\Hat{F}(\bx)-\II(\bx): \bx\in\sX\}$
for all $\Hat{F}\in\sC_b(\sF)$. It remains to observe that since $\sF$ is a \v{C}ech-Stone compactification,
every function $F\in\sC_b(\sX)$ is a restriction to $\sX$ of some $\Hat{F}\in\sC_b(\sF)$,
see \cite[IV.6.22]{Dunford-Schwartz(1958)}.
Therefore \rf{(1.7)} holds true for all $F\in\sC_b(\sX)$.

To prove that the rate function is tight,  suppose that there is $a>0$ such that
$\II^{-1}[0,a]$ is not compact. Then there is $\delta>0$ and a sequence $\bx_n\in\sX$
such that $d(\bx_m,\bx_n)>\delta$ for all $m\ne n$. Since Polish spaces have Lindel\"of property,
there is a countable number of open balls of radius $\delta/2$ which cover $\sX$. 
For $k=1,2,\dots$, denote by  $B_k\ni\bx_k$ one of the balls that contain $\bx_k$,
 and let $\phi_k$ be a bounded
continuous function such that $\phi_k(\bx_k)=2a$ and $\phi_k =0$ on the
complement of $B_k$. Then $F_n=\max_{k\geq n} \phi_k\searrow 0$ pointwise. On the other hand
\rf{(1.7)} implies
$\LL(F_n)\geq L_0+F_n(\bx_n)-\II(\bx_n)\geq L_0+a$,  contradicting
\rf{V1.4}.
\end{proofof}
\begin{Lemma} If $\LL(\cdot)$ is a \LVF\ then
\begin{equation}\label{(2.2)}
\inf_{\bx\in\sX}\{F(\bx)-G(\bx)\}\leq \LL(F)-\LL(G)
\end{equation}
\end{Lemma}
\begin{proof}
Let 
$\const=\inf_\bx \{F(\bx)-G(\bx)\}$ Clearly, 
 $F\geq G+ \const$. By positivity condition
\rf{V1.2} this implies  $\LL(F)\geq \LL(G+\const)=\LL(G)+\const$.
\end{proof}
The next lemma is implicitly contained in the proof of
\cite[Theorem T.1.1]{Bryc90}.
Let $\calP_{a} (\sX)$ denote all regular finitely-additive probability measures on
$\sX$ with the Borel field.
\begin{Lemma}\label{Lemma-T.1.1}
\label{T.1.1} If  $\LL(\cdot)$ is a convex \LVF\ on
$\sC_b(\sX)$, then there exist a lower semicontinuous function $\JJ:
\calP_{a} (\sX)\to[0, \infty ] $ such that
\begin{equation}\label{(1.4)}   \LL (F)=\LL(0)+\sup
\{ \mu(F)-\JJ(\mu): \mu\in
\calP_{a} (\sX)\}, 
\end{equation} and the supremum is attained. 

\end{Lemma}
\begin{proof} Let $\JJ(\cdot)$ be defined by
\begin{equation}\label{(1.5)}    \JJ(\mu)=\LL(0)+\sup \{ \mu(F)-\LL (F): F\in \sC_b(\sX)\} .
\end{equation}
and fix $F_0\in \sC_b(\sX)$. Recall that throughout this proof we assume $\LL(0)=0$. 

By the definition of $\JJ(\cdot)$, we
need to show that 
\begin{equation}\label{(2.7)}   \LL(F_o)= \sup_\mu \inf_F \{ \mu(F_0)
-\mu(F)+\LL(F)\},\end{equation} where the supremum is taken over all $\mu\in
\calP_{a} (\sX)$ and the infimum is taken over all $F\in \sC_b(\sX)$.
 Moreover, since \rf{(1.5)} implies that $\JJ(\mu)\geq \mu(F_0)-\LL(F_0)$
for all $\mu\in \calP_{a} (\sX)$, therefore $\LL(F_0)\geq \sup_\mu
\inf_F \{ \mu(F_0) -\mu(F)+\LL(F)\} $. Hence to prove \rf{(2.7)}, it
remains to show that there is $\nu\in \calP_{a} (\sX)$ such that 
\begin{equation}\label{(2.8)}   \LL(F_0) \leq \nu(F_0) -\nu(F)+\LL(F) \hbox{ for all $F\in \sC_b(\sX)$.}
\end{equation} 
(also, for this $\nu$, the supremum in \rf{(1.4)} will be attained) To find $\nu$, 
consider the following sets. 
Let   $$\calM= \{ F\in \sC_b(\sX): \inf_\bx [F(\bx)-F_0(\bx)]>0\}$$
and let $\calN$ be a set of all finite convex combinations of functions
$g(\bx)$ of the form $g(\bx)= F(\bx)+ \LL(F_0)- \LL(F)$, where $F\in
\sC_b(\sX)$.

It is easily seen from the definitions that $\calM$ and $\calN$ are
convex; also $\calM\subset\sC_b(\sX)$ is 
 non-empty since $1+F_0\in \calM$, and open
since $ \{F: \inf_\bx [F(\bx)-F_0(\bx)]\leq0\}\subset\sC_b(\sX)$ is closed.
Furthermore, $\calM$ and $\calN$ are disjoint. Indeed, take arbitrary
   $$
\calN\ni g= \sum \alpha_k F_k + \LL(F_0) - \sum \alpha_k\LL (F_k) . $$ 
Then   $$\inf_x\{ g(\bx)-F_0(\bx)\}
=$$
   $$\inf_x\{ \sum \alpha_kF_k(\bx) -F_0(\bx)\} - \sum \alpha_k\LL (F_k)+
\LL(F_0) \leq $$   $$\inf_x\{ \sum \alpha_kF_k(\bx) -F_0(\bx)\} -\LL(
\sum
\alpha_kF_k)+
\LL(F_0) \leq 0,$$ where the first inequality follows from the convexity of $\LL(\cdot)$ and
the
second one follows from \rf{(2.2)} applied to $F= \sum \alpha_kF_k(\bx)$ and $G=
F_0$. 

Therefore $\calM$ and $\calN$ can be separated, i.~e. there is a non-zero 
linear functional $f^*\in \sC_b^*(\sX)$ such that for some
$\alpha\in\sR $
\begin{equation}\label{(2.9)}   f^*(\calN) \leq \alpha< f^*(\calM), \end{equation}
see e.~g.  \cite[V.
2. 8]{Dunford-Schwartz(1958)} 

{\bf Claim:} $f^*$ is non-negative. 

Indeed, it is easily seen that $F_0(\cdot)$ belongs to $\calN$, and, as a
limit of $\eps+ F_0(\bx)$ as $\eps\to 0$, $F_0$ is
also
in the closure of $\calM$. Therefore by \rf{(2.9)} we have
$\alpha=f^*(F_0)$. To end the proof take arbitrary $F$ with $\inf_\bx
F(\bx)>0$.
Then
$ F+ F_0\in \calM$ and by \rf{(2.9)}
$$f^*(F)= f^*(F+F_0)- f^*(F_0)> \alpha- f^*(F_0)= 0.$$
This ends the proof of the claim. 

Without loosing generality, we may assume $f^*(1)= 1$; then it is well
known, see e.~g. \cite[Ch. 2 Section 4 Theorem 1]{Bergstrom(1982)}, that
$f^*(F)= \nu(F)$ for some $\nu\in \calP_{a} (\sX)$;
for
regularity of $\nu$ consult
\cite[IV.6.2 ]{Dunford-Schwartz(1958)}.
 It remains to check that $\nu$ satisfies \rf{(2.8)}. To this end observe
that since $F+\LL(F_0)- \LL(F)\in \calN$, by \rf{(2.9)} we
have
$\nu(F)+\LL(F_0)-
\LL(F) \leq \alpha= \nu(F_0)$ for all $ F\in \sC_b(\sX)$. This ends the proof
of \rf{(1.4)}.
  \end{proof}
\begin{proofof}{Theorem \protect{\ref{mini-T.1.1}}}
 Lemma \ref{Lemma-T.1.1} gives the variational representation \rf{(1.4)} with the supremum
taken over a too large set. To end the proof
we will show that $\JJ(\mu)=\infty$  on
measures $\mu$ that fail to be countably-additive. 

Suppose that $\mu$ is additive
but
not countably additive. Then Daniell-Stone theorem implies that
there is $\delta>0$ and a  sequence $F_n\searrow 0$ of
 bounded continuous functions
such that $\int F_n d\mu>\delta>0$ for all $n$. By \rf{(1.5)} and 
 $\sigma$-continuity 
$\JJ(\mu)\geq \LL(0)+C\int F_n d\mu -\LL(C F_n)\geq \LL(0)+ C\delta -\LL(C F_n)
\to \LL(0)+ C\delta$. Since $C>0$ is arbitrary, therefore
 $\JJ(\mu)=\infty$ for all $\mu$ that are additive
but not countably-additive. Thus \rf{(1.4)} implies \rf{CVX}.
\end{proofof}


\begin{thebibliography}{99}
\bibitem{Akian99} M. Akian,  Densities of idempotent measures and large deviations. 
Trans. Amer. Math. Soc.
351 (1999), 4515--4543.
\bibitem{Bergstrom(1982)} H. Bergstr\"om, {\sl Weak convergence of measures}. Acad. Press, New York, 1982. 
\bibitem{Bryc90} W. Bryc
 { On the large deviation principle by the asymptotic value method}.
 In: {\sl Diffusion Processes and Related Problems in Analysis}, Vol. I, ed.
M. Pinsky, Birkh\"auser, Boston 1990, 447--472.

\bibitem{deAcosta(1985)} A. de Acosta, Upper bounds for Large Deviations of
Dependent Random Vectors. Zeitsch. Wahrscheinlichk. Theor. Verw. Gebiete 69 (1985),  551--565.

\bibitem{D-Z} A. Dembo \& O. Zeitouni, {\sl Large Deviations
Techniques and Applications}.
Jones and Bartlett, Boston,  1993.

\bibitem{D-S} J-D. Deuschel \& D. W. Stroock, {\sl Large Deviations}.
 Pure and Applied Math vol. 137,
Academic Press, Boston, 1989.

\bibitem{Dunford-Schwartz(1958)} N. Dunford \& J. T. Schwartz, {\sl Linear Operators I}.
Interscience, New York, 1958.  

\bibitem{Dupuis-Ellis(1997)} P. Dupuis \& R. S. Ellis, {\sl A Weak
Convergence Approach to the Theory of Large Deviations}. Wiley, New York, 1997. 

\bibitem{OBrien(1996)}  G. L. O'Brien, Sequences of capacities, with connections to large-deviation theory. 
J. Theoret. Probab. 9 (1996), 19--35.

\bibitem{OBrien-Vervaat(1995)} G. O'Brien \& W. Vervaat,
 Compactness in the theory of large deviations.
Stoch. Processes Appl. 57 (1995), 1--10.

 \bibitem{Puhalskii(1997)} A. Puhalskii,
Large deviations of Semimartingales: a Maxingale Problem Approach I. Stochastics 61 (1997),
141--243. 

\bibitem{Varadhan(1966)} S. R. S. Varadhan, Asymptotic probabilities and differential
equations. Comm. Pure Appl. Math. 19 (1966), 261--286. 

\end{thebibliography}
\end{document}